\documentclass[12pt]{article}
\usepackage{amsmath}
\usepackage{amsfonts}
\usepackage{amssymb}
\usepackage{graphicx}
\usepackage{comment}
\usepackage{float}
\usepackage{textcomp}

\gdef\be{\begin{equation}}
\gdef\ee{\end{equation}}

\begin{document}
\setcounter{page}{0} \topmargin 0pt
\renewcommand{\thefootnote}{\arabic{footnote}}
\newpage
\setcounter{page}{0}


\begin{center}
{\Large {\bf Conformality in Hawksmoor's Ceiling,\\ 
\vspace{5mm}
Mercator's Projection and the Roman Pantheon}}
\vspace{1cm}

{\large John Cardy$^*$\\}  \vspace{0.5cm} 
{All Souls College, University of Oxford,\\ Oxford OX1 4AL, United Kingdom\\}
\vspace{2cm}
August 2026

\end{center}

\vspace{1cm}

\begin{abstract}
The dome of the Roman Pantheon is coffered with ribs surrounding sunken \em lacunaria\em, thus forming a grid. How this is achieved given the curvature of the dome has long been a subject for study and speculation. Although  detailed measurements now exist, thus far no single principle has emerged which fixes the overall geometry. Similar coffering occurs in Hawksmoor's design for the ceiling of the Buttery in All Souls College, Oxford. Both these examples are doubly curved, making their tiling with (almost) square coffers problematic. We address this using methods of differential geometry. Observing that the ribs intersect orthogonally, we hypothesize that they are images under a conformal (angle preserving) mapping of a uniform $M\times N$ tiling of a planar rectangle with squares. This is unique and is the inverse of  Mercator's projection of the ceiling to the rectangle. It minimizes an energy functional which captures the sum of the elastic and gravitational energies in the long wavelength limit. Thus, \em any \em construction method which is stable to long wavelength deformations necessarily leads to a conformal coffering. This then gives quantitative predictions with no adjustable parameters for the relative sizes and locations of each coffer,  which are in good agreement with photographic data, and consistent with direct measurements of the Pantheon. We suggest a protocol by which Hawksmoor's ceiling might have been constructed without advanced mathematics.  
\end{abstract}

\section{Introduction.}\label{secintro}
\subsection{Hawksmoor's Ceiling.}\label{sechawk}

Nicholas Hawksmoor (1661--1736) began his architectural career as a clerk under Sir Christopher Wren and later became a leading figure of the English Baroque movement. Although he is perhaps best known for his series of London churches, he
also conceived plans for extensive building works in Oxford and elsewhere.

In All Souls College he rebuilt the whole North end of the College, including the Hall and its adjacent Buttery. 
Although the exteriors were designed in a somewhat flamboyant neo-Gothic style to match the rest of the North Quadrangle, the interiors are decorated in a more restrained English Baroque manner. In earlier times, the Buttery would have been a place to store food and drink (in `butts'.) Until quite recently, however, it was where Fellows would take lunch, in somewhat crowded but convivial surroundings.

\begin{figure}[h]
\centering\includegraphics[width=0.41\textwidth]{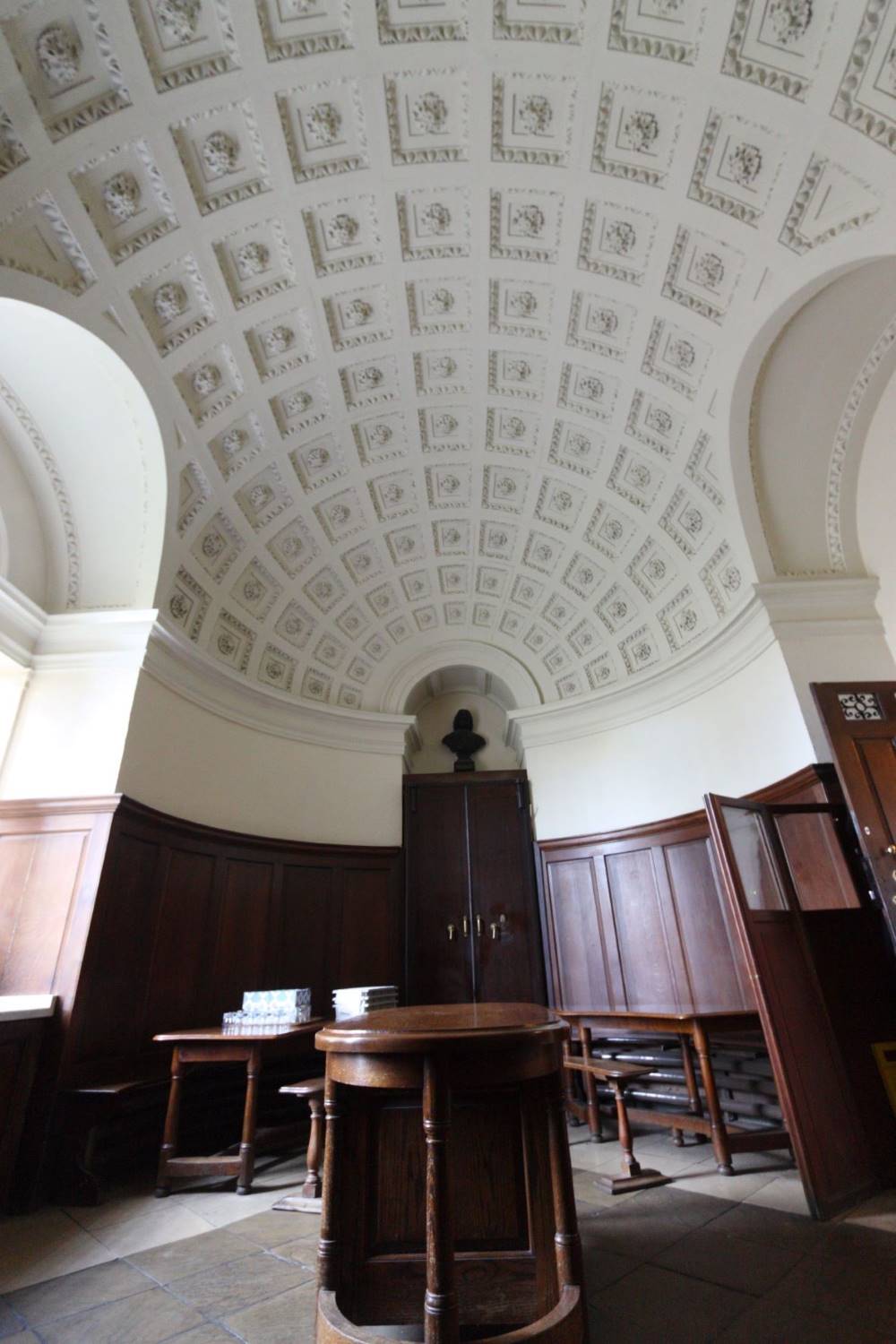}\hspace{2mm}
\includegraphics[width=0.5\textwidth]{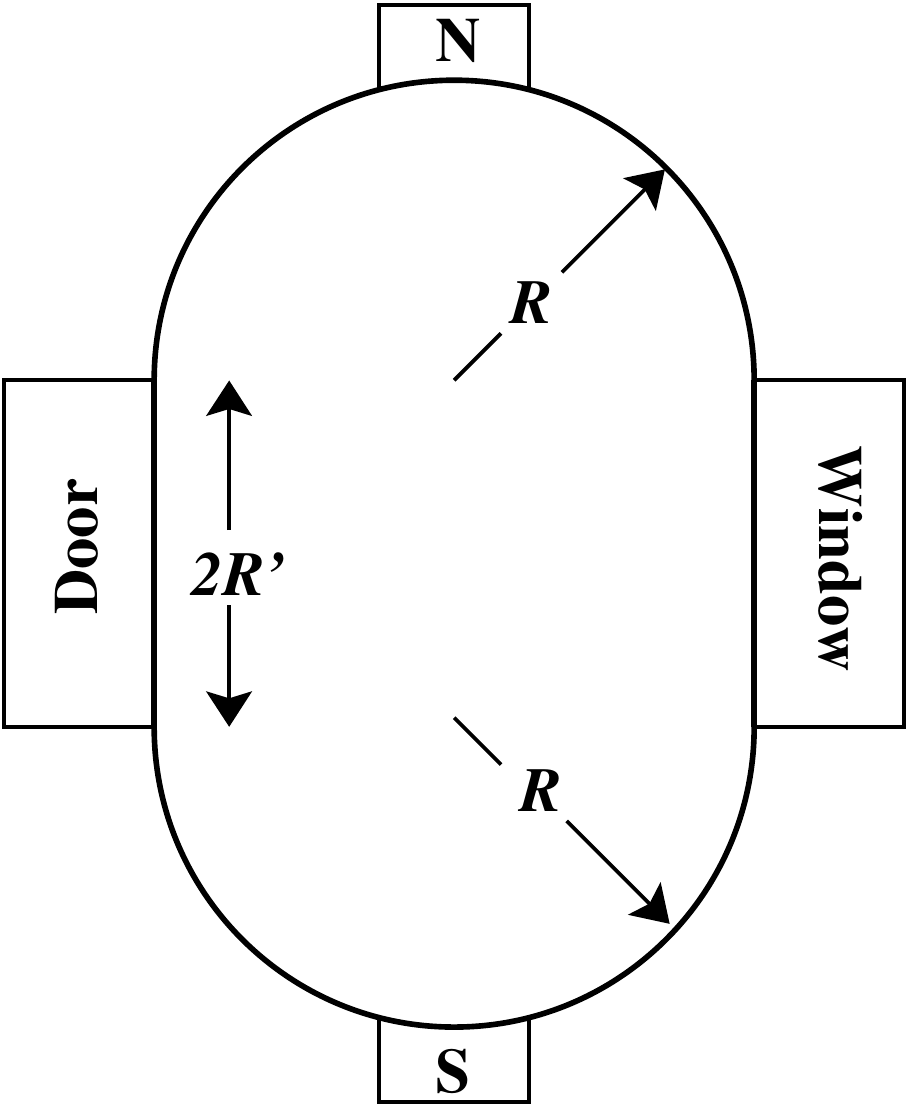}
\caption{\label{fig1}
Left: interior of the Buttery looking South. The oval plan of the walls and floor is apparent, as is the domed nature of the ceiling. The window and doorway intrude on this structure. Hawksmoor's coffered design for the ceiling uses approximately square coffers which vary in size and orientation, but whose edges always meet at right angles. The size of the grid diminishes at the polar ends, where it is interrupted by arched niches (Photo All Souls College).\\
Right: simplified plan of the Buttery at the level of the base of the ceiling. It consists of two half-discs connected by a rectangular central section. The ceiling surface is generated by rotating (say) the left half of the perimeter of the base through 180$^\circ$ about the axis of symmetry. The proportions are true: the coffering of the ceiling constrains the ratio $R'/R$ to be a rational multiple of $\pi$. In fact it is $2\pi/11\approx 0.57$.}
\end{figure}

\begin{figure}
\centering
\includegraphics[width=0.9\textwidth]{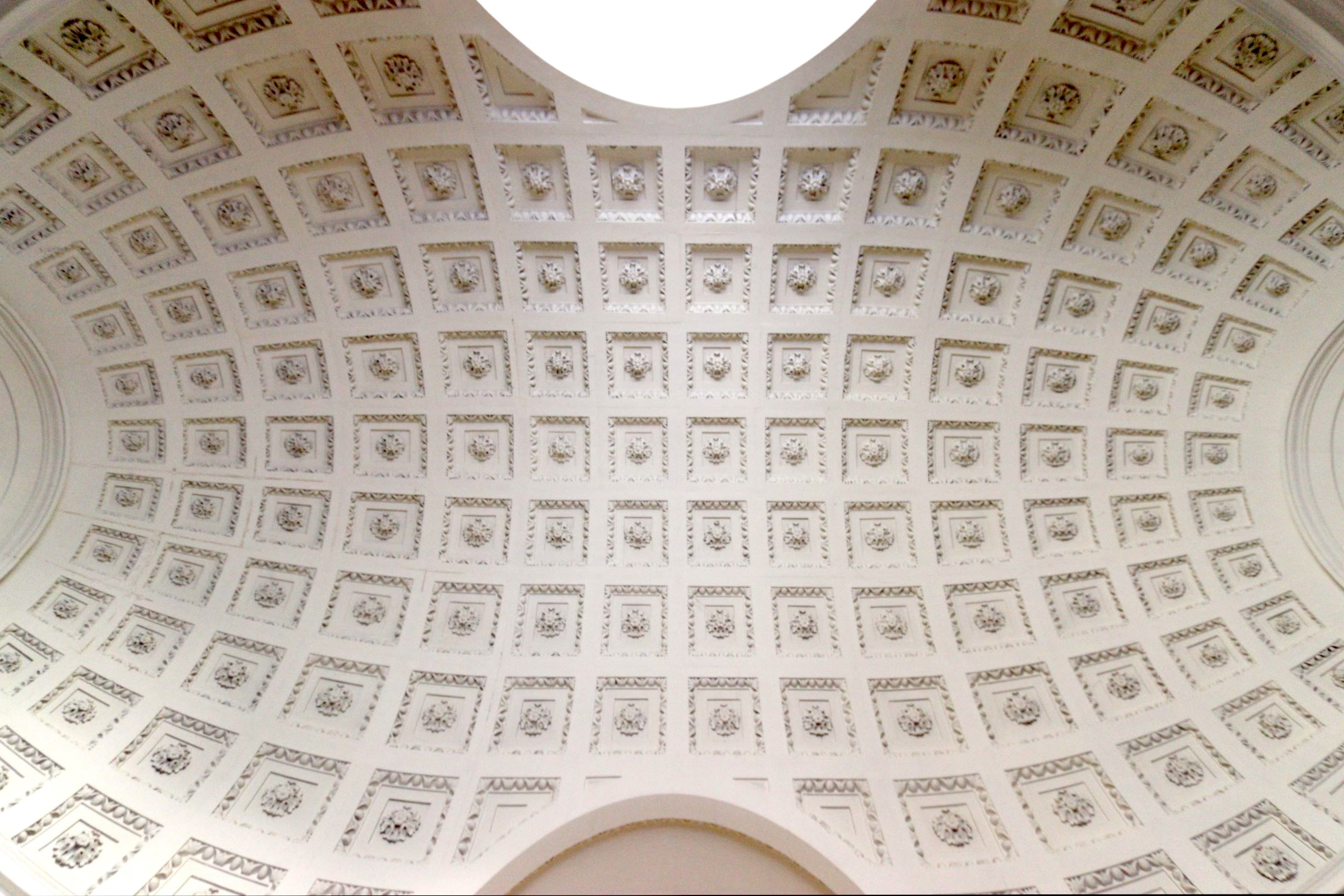}
\caption{\label{fig4} Inferior view of the ceiling. The image is oriented so that the doorway is at the top. It is distorted by the camera, so the orthogonality of the crossings of the coffer lines is not apparent except near the apex (Photo All Souls College.)} 
\end{figure}

The general layout of the Buttery and its ceiling is apparent from Figs.~\ref{fig1} and \ref{fig4}. External considerations constrained the room to have an oval shape, with its major axis oriented roughly North-South. It is interrupted on the East and West sides by a window and the entrance doorway respectively. The walls meet the ceiling at a horizontal cornice: thus its base is an oval, save where the window and doorway intrude in arches, and at the polar ends, which are terminated by smaller arches, each of which forms the roof of a niche.  Since the E-W arches clearly interrupt the pattern of coffers, it is useful to imagine them as absent, with the ceiling continuing down to meet the extended cornice. Although a purely elliptical oval base might have been more aesthetic, as in some of Borromini's Italian works, it appears in fact to consist of two semicircles at the polar ends, attached to straight sections on either side, as in a stadium. The ceiling itself is thus a half-capsule: two quarter spheres attached to a half-cylindrical barrel vault. At least this geometry serves to simplify some of the mathematics.

In the Baroque style it was customary to decorate such ceilings with coffers, that is
beams or arches forming a regular grid, surrounding sunken lacunaria in the shape of squares or hexagons. For definiteness we shall refer to each lacunarium, including all the structure up to the midlines of its bounding ribs, as a single coffer. The midlines then form smooth curves across the whole ceiling, which we call coffer \em lines\em.  In Hawksmoor's design these intersect orthogonally at the corners of adjacent coffers, which will be termed \em vertices\em. Our analysis will focus on the geometry of the lines and vertices rather than the structure within each coffer. Thus the problem becomes one of \em tiling \em a curved surface.

It is useful at this point to recall a few basic properties of curved surfaces smoothly immersed in three dimensional space. Near any given point on the surface it may be approximated by an ellipsoid, part of a sphere which has been uniformly stretched in one direction tangent to the surface and compressed in the other perpendicular direction. A transverse section of the surface in a plane parallel to each of these two principal directions may then be approximated by a circle. If these circles have radii $(R_1,R_2)$, the principal curvatures are their inverses  $\kappa_1=1/R_1$
and $\kappa_2=1/R_2$. The local curvature of the surface is then conventionally expressed in terms of the mean curvature $H=\frac12(\kappa_1+\kappa_2)$ and the Gaussian curvature $K=\kappa_1\kappa_2$. Note that if $K$ is non-zero then both $\kappa_1$ and $\kappa_2$ are non-zero. Such a surface is said to be doubly curved.

\begin{figure}
\centering\includegraphics[width=0.6\textwidth]{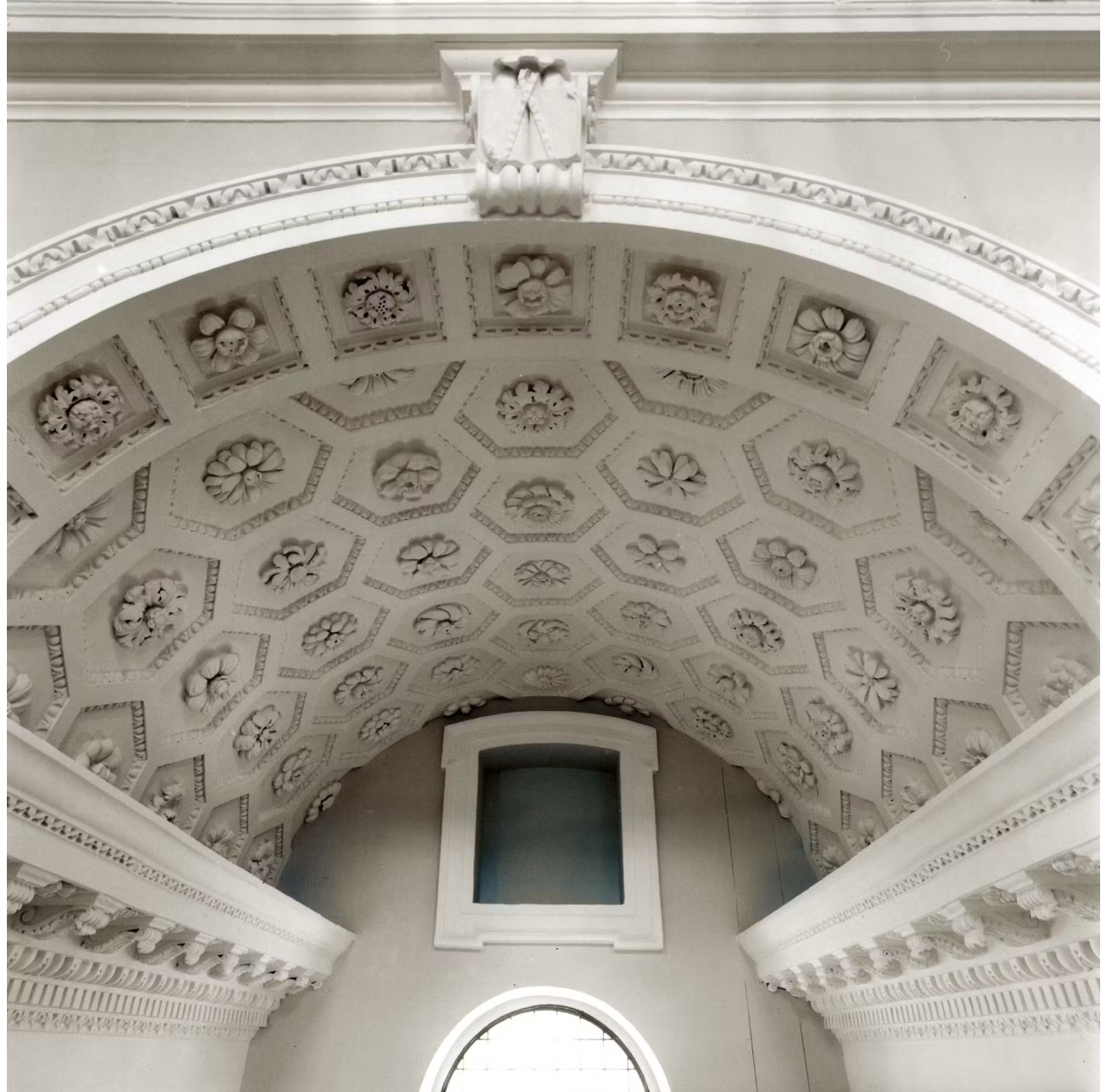}
\caption{\label{figS} Hawksmoor's barrel vaulted coffered ceiling in Christ Church at Spitalfields. Since the surface has no Gaussian curvature it may be tiled uniformly with regular hexagons or squares.}
\end{figure}

 Coffering initially served to support the structure while minimizing weight, but later it became purely decorative. Examples are the arches lining the nave of Hawksmoor's masterpiece, Christ Church at Spitalfields (Fig.~\ref{figS}). These have a cylindrical barrel form, and therefore no intrinsic, Gaussian, curvature: that is, one may imagine unrolling them to be flat, without disturbing the relative positions of the coffers. Since a uniform tessellation of a plane with squares or hexagons is always possible, it follows that this holds also for a barrel vault. 

However, the Buttery ceiling has non-zero Gaussian curvature: it is doubly curved, and cannot be flattened out in the same way.\footnote{It is an intrinsic property because a geometer living within the surface, on measuring the ratio of the circumference to the radius of a small disc embedded therein, would find it unequal to $2\pi$, and would conclude that their geometry is non-Euclidean.} It is not, in general, possible to cover such a surface with identical squares without leaving some gaps between them.  Most coffered domes meet this requirement by allowing trapezoidal coffers, although mathematically it is always possible to tessellate almost all of the surface with almost square coffers as long as they are allowed to vary in size.

For whatever reasons: aesthetic, in homage to classicism, or on strictly practical grounds, Hawksmoor chose to relax the condition that all coffers be identical in size, while maintaining their square shape to the greatest extent possible. We shall call such a design a \em Hawksmoor coffering\em.\footnote{Although there seems to be only one other example of such a doubly curved ceiling with designed by him, in the Mausoleum at Castle Howard, which has a dome similar to that of the Pantheon.}
His problem then consisted in determining  how their size and position should depend on the overall geometry of the ceiling. There are two levels at which to pose this: the first is the precise theoretical determination of the ideal relative coordinates of each coffer in three dimensional space, and the second is the more practical question of devising a constructive scheme which achieves a sufficiently faithful approximation to this, preferably without any advanced mathematics. In this note we mainly focus on the first, giving an explicit answer for any shaped ceiling which is a surface of revolution. However, the argument points the way to a possible solution to the practical problem.

\subsection{The Roman Pantheon.}\label{secpanth}

\begin{figure}
\centering
\includegraphics[width=0.41\textwidth]{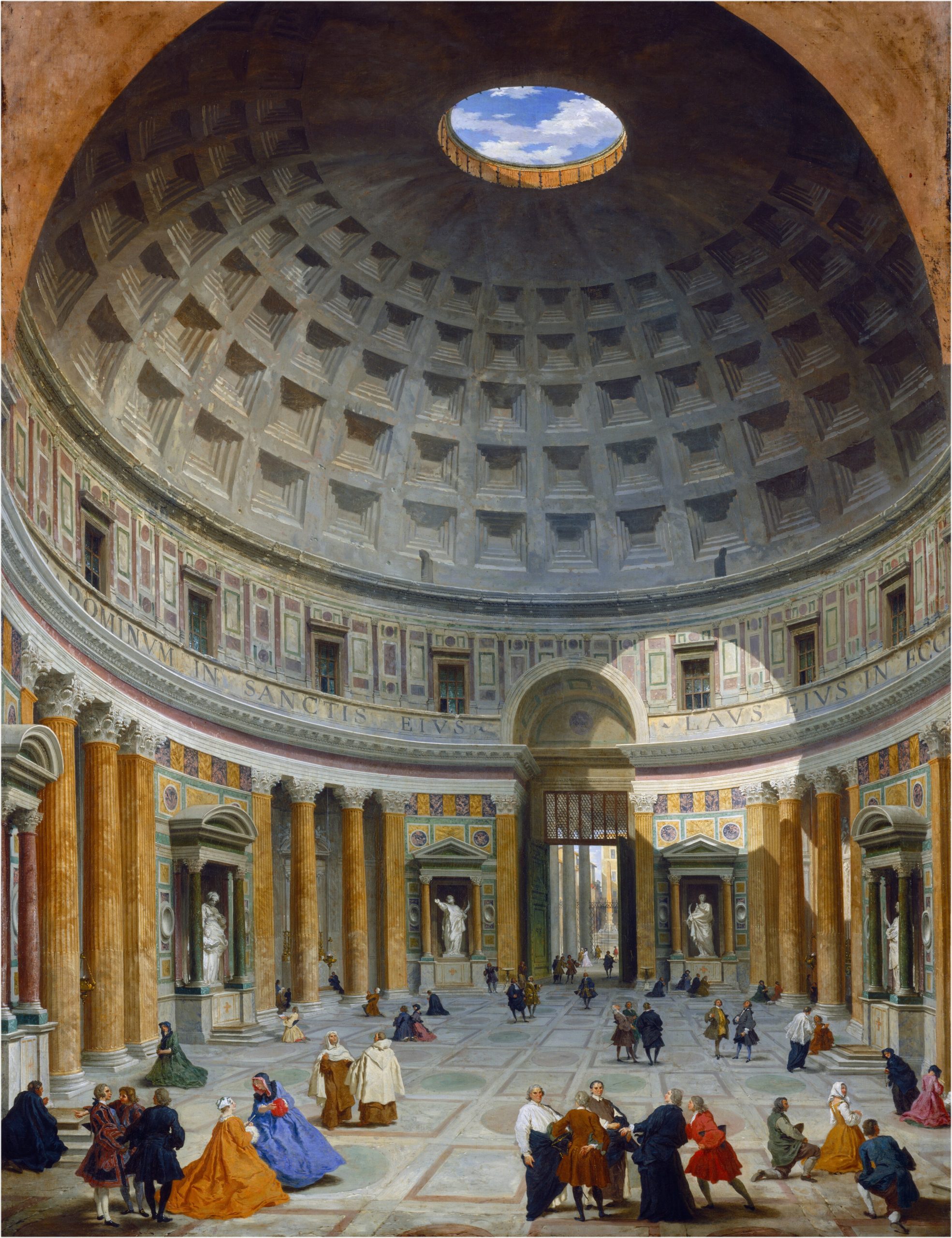}
\hspace{2mm}
\includegraphics[width=0.52\textwidth]{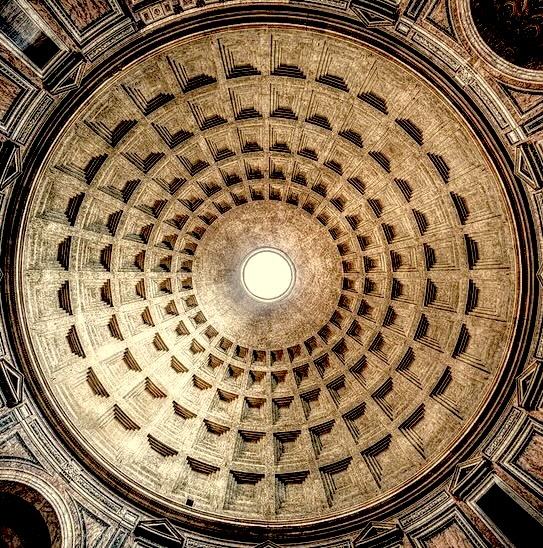}
\caption{Left: Giovanni Paolo Panini, \em Interior  of the Pantheon (c.~1734) \em (National Gallery of Art).\\
Right: \label{panphoto} Modern inferior view of the Pantheon dome. (\textcopyright\ David Lown)}
\end{figure}

One of the earliest examples of the square coffering of a dome is in the Pantheon in Rome (125 CE), shown in Fig.~\ref{panphoto}. Hawksmoor was fascinated by classical forms, so there is every reason to believe he had studied drawings of this, although apparently he never travelled to Italy.  In fact a comparison of Fig.~\ref{fig1} and Fig.~\ref{panphoto}suggests that bisecting the Pantheon dome to give a quarter-sphere, which is then rotated through 90$^\circ$ so that its vertical axis lies in a horizontal plane, gives a pattern of coffers which strongly resembles that of the quarter-spheres at each end of the Buttery ceiling, with the role of the oculus played by the N-S niches. One difference, which, however, turns out to be immaterial, is in the direction of the gravitational force relative to the axis of symmetry.  It should also be noted that, while the Buttery has a rather small ceiling, and Hawksmoor's coffering is quite simple, the Pantheon dome consists of several layers which, apparently, serve to support the whole massive structure composed of concrete of varying thickness and composition. Thus one should not expect it to conform to any simple rules in the same way, although  in fact it appears to do so. 
 
Unlike the All Souls Buttery ceiling, the geometry of the Pantheon dome has been the subject of numerous studies over the years (see Sec.~\ref{secrec} for a summary of recent literature), and has inspired many copies, both in the Roman era\cite{Aun2025} and in modern times,  notably Thomas Jefferson's rotunda at the University of Virginia, and his memorial in Washington, DC\cite{Fle2003}.

Most investigations of the Pantheon focus on its aesthetic, historical, and numerological significance, often using rather inaccurate measurements based on artists' impressions. More recently more reliable information has become available using modern metrological methods, which, however, tends to point out the deviations from general laws  rather than their underlying universality. The approach described  here emphasises the latter, yielding explicit predictions for the mean sizes and positions of the coffers rather than their variability.

\subsection{Mercator's Projection}\label{secmerc}

Hawksmoor's problem of constructing a square grid on an intrinsically curved surface has an inverse, that of mapping a  curved surface onto a plane. 
This had arisen some time earlier in the age of global exploration, with the need to make a planar chart which would faithfully represent a sufficiently large portion of the Earth's surface. It is not possible to preserve relative distances between different points while maintaining the angular relationships required to plot a course between them. Since the latter was more important for compass navigation, Gerardus Mercator's solution to this problem in 1569 using his famous projection  became a standard method for producing sailing charts, and is still used today for GPS based navigation maps, although as a representation of the entire globe it lends too much emphasis to the polar regions. Since it too locally preserves angles, it is also an example of a conformal mapping. Since it maps a portion of the globe to a rectangle whose sides are lines of longitude and latitude, and right-angled corners map to the same, it is unique and therefore must be the inverse of the Hawksmoor mapping to a portion of a spherical globe. Similarly Mercator's projections of the whole Earth (excluding the polar regions) are to a finite cylinder, and are the inverse of the conformal mapping which gives the coffering of an extended dome. 
This holds even for a non-spherical globe, as long as there is axial symmetry. Just as  Mercator's projection amplifies features near the poles like Greenland, its inverse has the opposite effect, consistent with Hawksmoor's squares diminishing in size towards the N-S ends of the ceiling.

\begin{figure}
\centering
\includegraphics[width=0.4\textwidth]{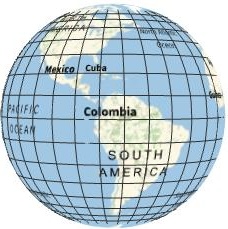}
\includegraphics[width=0.15\textwidth]{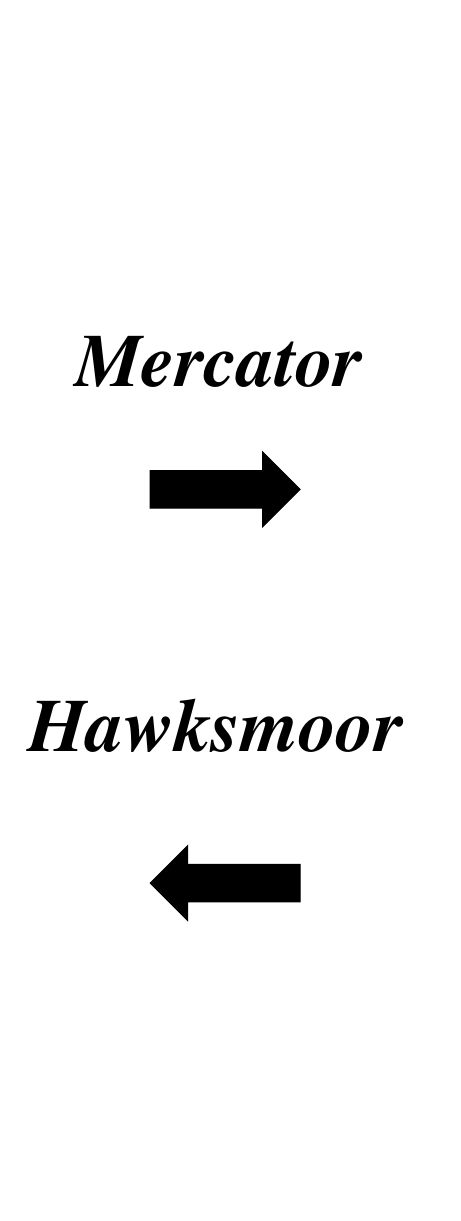}
\includegraphics[width=0.4\textwidth]{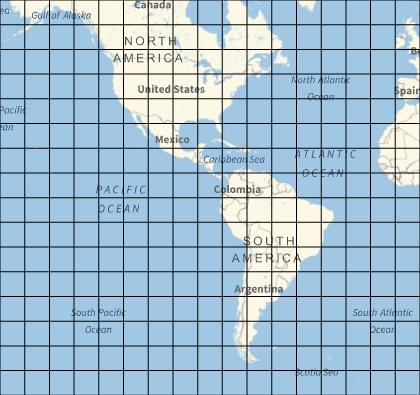}
\caption{\label{figmerc}Hawksmoor's coffering is the inverse of Mercator's projection: if a portion of the Earth's globe is superimposed onto a Hawksmoor coffered hemispherical ceiling, under Mercator's projection the continents are mapped as usual, and the coffer lines map into a regular square grid. Note that the parallels are not spaced in equal steps of the latitude angle $\theta$, as for the standard illustrations, but rather to create a square tiling on the map, but they do faithfully represent the local compass directions (plumb lines).}
\end{figure}

Thus if one takes a part of Mercator's projection of the continents and rules it with a uniform square grid, under the conformal mapping to a hemispherical ceiling the grid maps to Hawksmoor's coffering, while the continents map to their correct relative positions on a hemisphere of the globe. This is illustrated in Fig.~\ref{figmerc}.

\section{Methodology and main results.}\label{secmain}

One purpose of this note is to suggest an answer to the more general question ``What is the optimal way of coffering an intrinsically curved 2-dimensional surface immersed in 3-dimensional Euclidean space?" This of course depends on subjective criteria, but the most remarkable common features of the All Souls ceiling and the Pantheon dome are that the coffer ribs, where they cross,  do so orthogonally, and that they surround lacunaria which are close to being square, although their relative sizes differ. We shall assume, therefore, that achieving this was Hawksmoor's guiding principle, although it is by no means obvious that this is even possible in general. This point of view has recently been advocated by Fern\'andez-Cabo \cite{Fer2013}. We defer until Sec.~\ref{secDir} the more fundamental question of why this principle should apply. 

\begin{figure}
\centering
\includegraphics[width=0.8\textwidth]{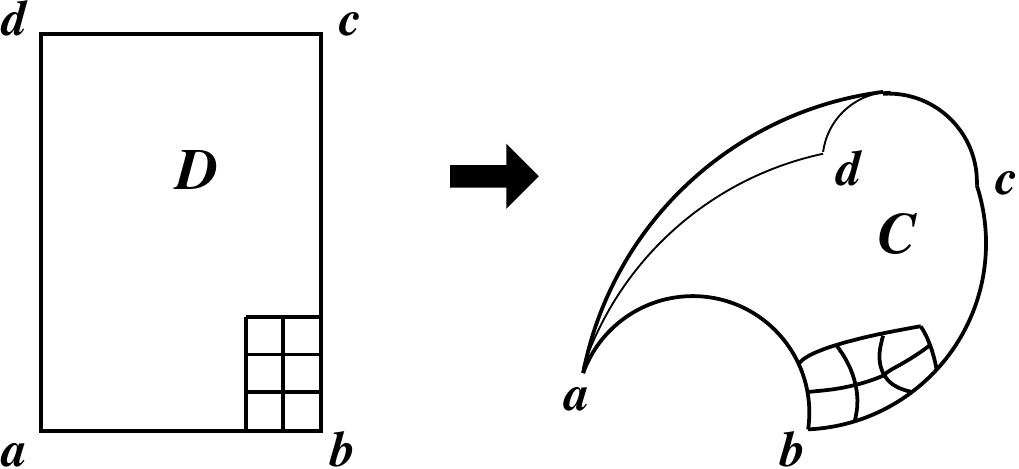}
\caption{\label{mapping1}
A general smooth map $f$ from a planar rectangle $D$ to a vaulted ceiling $C$. A uniform tiling of the rectangle with squares (a few of which are shown) maps to a general coffering of the whole ceiling with curvilinear quadrangles. If $f$ is conformal, however, the coffer lines  always intersect orthogonally, and the coffers are almost square.}
\end{figure}

An important feature of the approach adopted in this paper is that we  begin with a continuum description of the ceiling $C$ to be coffered. Suppose $(x_1,x_2)$ are Cartesian coordinates on a plane ${\mathbb R}^2$, and let $f:(x_1,x_2)\to \vec X=(X_1(x_1,x_2),X_2(x_1,x_2),X_3(x_1,x_2))$ 
be a smooth map from a fixed domain $D$ of the plane onto $C$, immersed in 3-dimensional Euclidean space. See Fig.~\ref{mapping1}. The map should be single-valued and invertible, that is for every point $(x_1,x_2)$ in $D$ there is a unique point $\vec X(x_1,x_2)$ on the ceiling, and \em vice versa\em. The
isolines along which $x_1$ is constant, and those of constant $x_2$, then form a curvilinear coordinate system\footnote{For this to be possible, $C$ should also be simply connected, which the case for the Hawksmoor ceiling. For the Pantheon, the map is between a finite length cylinder with a flat metric and the dome, a hemisphere minus the oculus. Both of these are topologically equivalent to an annulus, whose modulus is determined by the ratio $M/N$.} on  $C$.

The map also also gives a coffering design of $C$ with curved quadrangles, where the  coffer lines are given by $x_1=m\ell$ and $x_2=n\ell$ and the vertices lie at $\vec X(m\ell,n\ell)$, where $\ell$ is the mesh size and $m$ and $n$ are integers such that  $(m\ell,n\ell)$ lies within $D$. By fixing the map and letting $\ell\to0$, one obtains a sequence of increasingly fine cofferings of $C$ with curvilinear quadrilaterals. 
We may then pose the question: \em\lq\lq How might one characterise those maps which generate a Hawksmoor coffering of $C$?''\em

A theorem from differential geometry \cite{Lee2009} then implies that, for any sufficiently smooth ceiling shape, such a coffering is always possible. This is due to the existence, under very general conditions, of an \em isothermal \em coordinate system\footnote{So named because the equation describing heat conduction simplifies in these coordinates.} on the ceiling $C$, related to Cartesian coordinates  in a subset $D$ of the plane by a \em conformal mapping. \em This has the property that it locally preserves relative angles: if two smooth curves in $D$ intersect at some point with some angle between their tangents there, their images in $C$ meet at the same angle, even though lengths may be distorted. Thus a regular tiling of $D$ with identical squares is mapped to a Hawksmoor coffering of its image in $C$. It covers the whole of $C$ if the connected region in $D$ is made up of identical squares. The \em conformal hypothesis\em, put forward here, then asserts that this is unique and yields the observed real coffering design.

The theorem cited above is a simple consequence of some standard results in the theory of two-dimensional surfaces immersed in three dimensions, largely due to Gauss \cite{Gau1828} in the early 19th century, thus postdating Hawksmoor by a century or more.  The justification for this revisionist approach is that it shows how, given the principles stated above, Hawksmoor's solution is unique, whether he discovered it by a relatively simple iterative process (as outlined in Sec.~\ref{sechow}), or by more advanced mathematics. 

\subsection{Explicit results}
\begin{figure}
\centering
\includegraphics[width=0.4\textwidth]{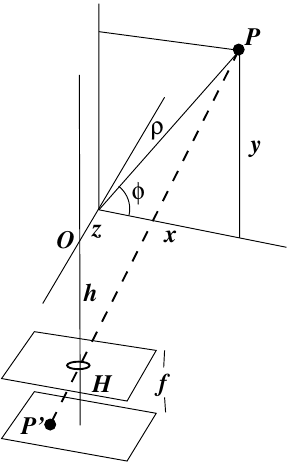}
\caption{\label{cylcoord}
Cylindrical coordinates used to specify the shape of the ceiling, with the $z$-axis horizontal in a N-S direction. Also shown is the
schematic setup for a pinhole camera (not to scale): a light ray from a point $P$ on the ceiling passes through a small hole $H$ in a horizontal screen situated a depth $h$ below the origin $O$, which is at the centre of the base of the ceiling. It is recorded at $P'$ where it hits a plate a distance $f$ below the screen.}
\end{figure}

For a generally shaped ceiling the identification of isothermal coordinates is a difficult but tractable problem \cite{Aul2013}. However it simplifies considerably if the ceiling is a surface of revolution, generated by rotating some simple curve in two dimensions about a fixed axis joining its two ends. This appears to be the case for the Buttery ceiling. Given a particular such shape, defined in cylindrical polar coordinates $(\rho,\phi,z)$ (see Fig.~\ref{cylcoord}) by a general equation $\rho=\rho(z)$, we are then able to obtain precise expressions for the locations of the vertices and coffer lines on the ceiling.

Our results for the ceiling problem are summarised below. More details are provided in Sec.~\ref{secmath}.
\begin{itemize}
\item For any sufficiently smooth shaped ceiling it is always possible to coffer it so that the coffer lines intersect orthogonally, and the coffers are approximately square, becoming exactly so in the limit when the grid becomes increasingly fine. 

\item However, such theoretical cofferings typically exhibit points at the boundary at which the size of the coffers vanishes. If one demands that the number of coffers is finite, so they have a minimum size, and that they fit without gaps along the boundary, the solutions fall into a discrete number of possibilities which are conformal images of a planar $M\times N$ rectangle tessellated with squares, of which the upper and lower edges map to the N-S arches, and the other edges to the boundary of the base, as illustrated in Fig.~(\ref{mapping1}).
 \end{itemize}
 
\begin{figure}
\centering
\hspace{-1.5cm}
\includegraphics[height=0.45\textwidth]{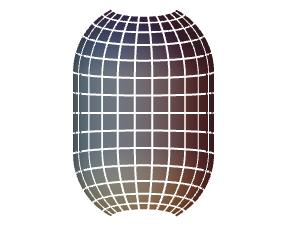}
\includegraphics[height=0.5\textwidth]{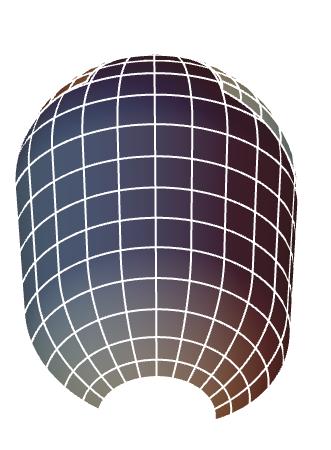}\\
\includegraphics[height=0.4\textwidth]{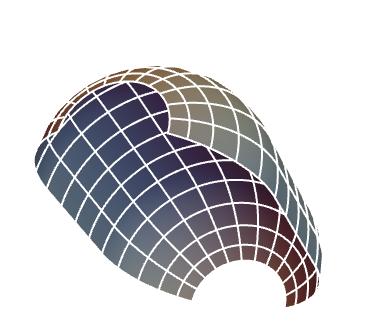}
\includegraphics[height=0.42\textwidth]{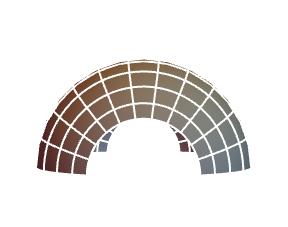}
\vspace{5mm}
\caption{\label{recon} Various 3D views of the coffering of the Buttery ceiling as predicted by the conformal hypothesis, using the actual number of coffers in each row and column. Only the coffer lines are shown. These are to be compared with Figs.~(1,2), taking into account the camera projection and the intrusion of the doorway and window arches.}
\end{figure}

If, in addition, the ceiling is a surface of revolution, \em i.e. \em its shape is independent of the azimuthal angle $\phi$ (see Fig.~\ref{cylcoord}), the solutions have the following properties:
\begin{itemize}
\item all the coffer lines are either in the N-S direction (meridians) or orthogonal to these (parallels.)
\item it follows that each E-W row of coffers between two adjacent parallels forms an arch, each with the same number $N$ of coffers, and that all those in a given arch are congruent, that is have the same size and shape. 
\item the mean width of any coffer in such an arch is then given by $\pi/N$ times its radius, \em i.e. \em the mean height of the apex above the base of the ceiling.
\end{itemize}
Thus the only remaining variables are the spacings between neighbouring arches. These are given by the following result:
\begin{itemize}
\item
For any surface of revolution, specified in cylindrical polar coordinates by $\rho=\rho(z)$, the coordinates of the $(m,n)$ vertex in an array with $N$ coffers per E-W row are $\phi=n\pi/N$, $z=z_m$ and $\rho=\rho(z_m)$, where 
\begin{equation}\label{zm}
m\pi/N= \int_0^{z_m}\frac{\sqrt{1+\rho'(z)^2}}{\rho(z)} \,dz\,.
\end{equation}
\item The size (height and mean width) of the coffers in the $m$th row is then proportional to $\rho(z_m)$.
\end{itemize}
For a hemispherical ceiling of radius $R$, this predicts the cylindrical coordinates of the $(m,n)$ vertex to be
\begin{equation}\label{coords}
z_m=R\tanh(\pi m/N)\,,\qquad \rho_m=R\,{\rm sech}(\pi m/N)\,,\qquad\phi=\pi n/N\,.
\end{equation}
In this special case, there is an iterative solution
\begin{equation}
z_{m+1}=\frac{z_1+z_m}{1+z_1z_m/R^2}\,.
\end{equation}
For the half-capsule shape of the Buttery ceiling, with the actual values of $M=15$ and $N=11$ used by Hawksmoor, the predicted coffering is shown in Fig.~\ref{recon}. 

\subsection{Comparison with observational data}
A second feature of our methodology is how we then compare these results with observations of the actual ceiling. Taking precise measurements of the ceiling itself runs up against the problem of defining exactly where the smooth background surface actually lies, as well as possible issues of accessibility, although remote sensing may provide a solution to the latter.
Instead, we have taken the path of comparison of the conformal predictions against mid-range photographic data, which, although less quantitative, has the advantage of being readily available and  being sensitive to features on a wide range of scales. The down side of this is that the camera projection distorts the coffer lines so that they no longer appear to meet orthogonally. This happens because the distance between the objective lens and the ceiling is
of the same order as the overall size of the ceiling. However this easily taken into account by applying a fixed projective transformation to the predicted position of each vertex.  Although the ratio between the above two length scales is usually not recorded, it appears in our calculations as an adjustable parameter and may be treated as such, chosen to obtain the best overall fit to the data.

\begin{figure}
\centering
\includegraphics[width=0.85\textwidth]{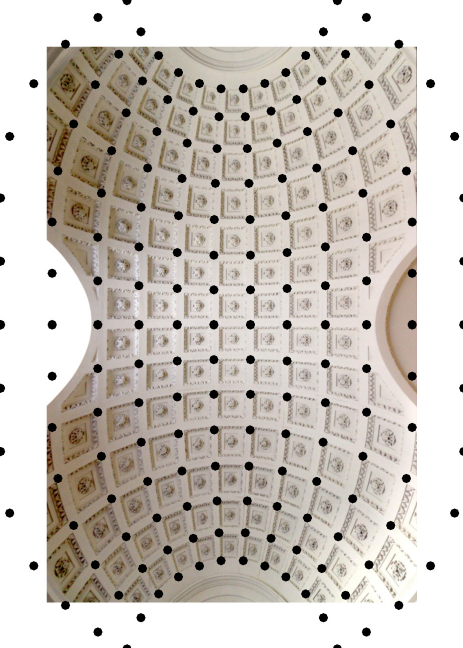}
\caption{\label{comp} The result of a theoretical camera projection of the predicted 3d locations of the vertices shown in Fig.~(\ref{recon}), overlaid onto the photographic image in Fig.~\ref{fig4}, expanded to show the fine details. Apart from the overall scale, there is one adjustable parameter, the height difference between the camera and the base of the ceiling, which was fixed to fit the curvature of the images of the N-S arches.  Only the coffer vertices are shown.}
\end{figure}

The results of this comparison between the results of the conformal hypothesis and the photograph in Fig.~\ref{fig4} are shown in Fig.~\ref{comp}. The agreement is excellent, the best fit being consistent with the camera being positioned approximately 150 cm above the floor.

The results for the dome are very similar, except that $z$ now labels the vertical axis and the mapping is from a finite cylinder, rather than a rectangle, to the surface of the dome. They may be compared, after applying the camera projection, to a photograph taken looking vertically upwards along the axis of symmetry, as in  Fig.~\ref{panphoto}. The result of comparing this with an actual photograph is shown in the lower figure in Fig.~\ref{comp}. Once again the agreement is satisfactory, although there are some systematic angular discrepancies, which could be attributed to the camera not being accurately aligned, or to  imperfect construction methods.

\begin{figure}
\centering
\includegraphics[width=0.98\textwidth]{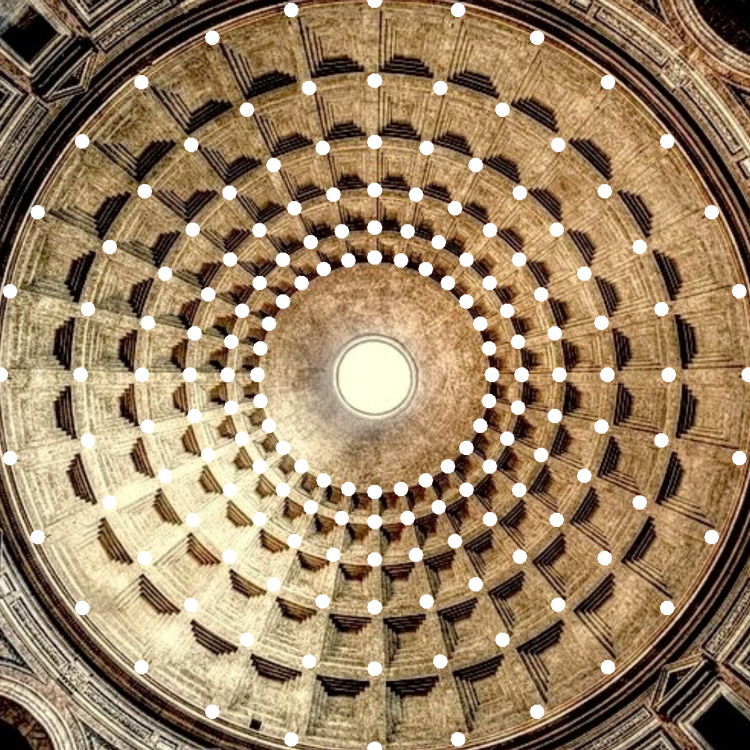}
\caption{\label{pancompph2} Expanded view of the right hand picture in Fig.~\ref{fig4}, showing
Pantheon dome with predicted positions in 3d of coffer vertices, as projected by the camera and then overlaid (white dots) onto a photograph (Fig.~\ref{panphoto}) of the actual dome. Again the ratio of height difference between the lens and the base of the dome to its radius has been adjusted to achieve the optimum overall fit. }
\end{figure}

\begin{figure}\centering
\includegraphics[width=.8\textwidth]{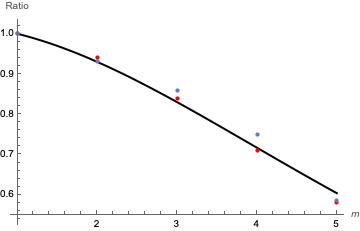}
\caption{\label{ratios} Pantheon: comparison between the measured values quoted in Table 2 of \cite{Ali2017} of the ratios $H_m/H_1$ (blue dots) and $W_m/W_1$ (red dots), where $H_m$ is the height and $W_m$ is the mean width of the $m$th row of coffers, and the conformal prediction (solid curve). Note that this is parameter-free: it is not a fit.}
\end{figure}

Until relatively recently the only direct comparison would have been against measurements extracted from photographs and artists' impressions. However, Ref.~\cite{Ali2017} reports  the results of modern metrological methods for the mean heights and widths of each row of coffers. As shown below, according to the conformal hypothesis, for a hemispherical ceiling these should be proportional to
$\text{sech}(2\pi m/N)$, where $N$ ($=28$ in this case) is the number of coffers in each row.\footnote{$\text{sech}(\xi)\equiv 2/({\text e}^\xi+{\text e}^{-\xi})$. Hyperbolic functions like this had already implicitly appeared in Mercator's work, although their connection to the number $e$ was understood only after the invention of calculus.}  Thus their ratios should be universal. A comparison against the measurements is shown in Fig.~\ref{ratios}. The agreement is fairly good, especially since the predicted curve has no free parameters. However, because the number of rows is limited, this does not probe the predicted exponential tail of the distribution. In addition, the irregularities in the coffering already observed in \cite{Ali2017} prevent the simple secant formula from being more than a rough fit.

\subsection{Comparison with other recent work.}\label{secrec}

To the author's knowledge, there has never, up to now, been a quantitative analysis
of the Hawksmoor ceiling, and unfortunately his original drawings are missing from the College's collection of his papers. A deep literature search of joint appearances of the terms ``coffering'' and/or ``Pantheon'' with ``conformal'' and/or ``Mercator'' yielded no results, suggesting that the proposal of using conformal mappings to generate coffering patterns, put forward in the present work, is new, and, likewise, the identification of these mappings with the inverses of Mercator-like projections.

The current scholarly monographs on the Pantheon and its construction are by Licht (1968) \cite{Licht1968} and Waddell (2008) \cite{Wad2008}. For a brief but excellent introduction see MacDonald (1976, repr. 1981, 2002) \cite{Mac1976}. A diverse compilation of the history \cite{Mar2015} was edited by Marder \& Wilson Jones (2015).

There is, of course, an extensive earlier literature on the design of the Roman Pantheon (see, for example, Palladio (1738) \cite{Pal1738}, Desgodets (1695) \cite{Des1695}, Rondolet (1860) \cite{Ron1860}.)
These authors have relied on contemporary measurements of varying degrees of reliability, which has not prevented them from embarking on flights of geometric, arithmetic,
musicological and astronomical fancy. (For a survey, see \cite{Fle2019}.) 
However in 2005 the Bern Digital Pantheon Project
(BDPP) was inaugurated, \cite{Alb2009} which, by 2019, had reportedly \cite{Fle2019} scanned the existing structure to a margin of error of 5 mm. These detailed measurements are useful in correcting previous assumptions and uncovering flaws in the original design, as well as the many improvements and repairs to the structure over
the millennia, but do not necessarily shed light on the basic principles used by the Romans.
Indeed Aliberti and Alonso-Rodr\'iguez \cite{Ali2017} state in their abstract
\em ``The configuration of the coffers is governed by a complex design that is based on geometric-constructive laws yet to be defined.'' \em  In their conclusions they state \em ``regarding the coffer shape, we propose the hypothesis that Roman constructors
designed the external form of the coffers with a square proportion based on the near
equivalence between the lowest edge and the vertical edge, adapting to the vertical
alignment of the external lines.'' \em Similarly Fern\'andez-Cabo \cite{Fer2013} approaches from the point of view of the designer \em``seeking  a geometric pattern of the
coffering so that it could be perceived by the human observer as having the shape most similar to a regular squared grid with each coffer appearing to be a perfect
square." \em This is close to the spirit of the present work. However, after reviewing the 
available data, he concludes that none of the existing theories correctly describe this, and 
\em ``the correct geometrical solution -- closest to the measurements actually constructed --
keeping this coffer regularity condition is unsolved in planar geometry of plans and
sections, as it would have to unfold the plan of the spherical segments; this is
geometrically and mathematically impossible due to the property of double curvature of
its surface. Single-curved surfaces can be unfolded but not double-curved surfaces. Thus,
the correct solution involves applying certain knowledge of spherical trigonometry.'' \em

We believe that the present work provides that solution, although we have chosen to express it in mathematical terms rather more advanced than those of spherical trigonometry. In fact, the conformal mapping methods we use are not new to the subject of architectural mathematics. (See, for example, the encylopaediac 2007 text by
 Pottman \em et al. \em\cite{Pot2007}.)

Other recent works (\cite{Fuc2023}, \cite{Aun2025}) have focused on the simple ratios of particular measured lengths. Since the plan of the whole edifice involves incorporating both spherical and cuboidal elements, these ratios should involve the number $\pi$, similar to the discussion in Sec.~\ref{secbutt} for the Buttery. However their integer arithmetic appears to work only if it is assumed that $\pi$ is taken to have the value $\frac{22}7$, and also seems to depend on the perimeter of the dome being divided into $N=28$ equal arcs. We emphasise that our analysis makes no such assumptions, and generates cofferings of the dome with any positive value of $N$, although the coffers become perfectly square only in the limit $N\to\infty$.

 \section{Mathematical Analysis} \label{secmath}
The non-mathematical reader may wish to omit this section. However it uses only well-known results in differential geometry \cite{Lee2009}, and the actual manipulations are elementary. We begin by recalling some basic facts about 2-dimensional manifolds.
\subsection{Differential geometry in two dimensions.}\label{secdiffgeom}
A 2-dimensional differentiable manifold is a set of points labelled by continuous real coordinates $(x_1,x_2)$, and endowed with a metric which specifies the squared distance $(ds)^2$ between infinitesimally separated points $(x_1,x_2)$ and $(x_1+dx_1,x_2+dx_2)$:
\begin{equation}\label{metric}
(ds)^2=g_{11}(dx_1)^2+g_{22}(dx_2)^2+2g_{12}(dx_1)(dx_2)\,,
\end{equation}
where the components $g_{ij}$ of the metric are in general functions of the coordinates. In general a coordinate system is not unique: we could use any suitably smooth functions $(x'_1,x'_2)$ of $x_1$ and $x_2$, as long as we also modify the metric in such a way that $(ds)^2$ is invariant. To be a manifold it should then locally look like the Euclidean plane: that is, for each point $P$ there exists a choice of coordinates  such that $g_{11}(P)=g_{22}(P)=1$ and $g_{12}(P)=0$. However in general such a choice simplifies the metric only exactly at $P$ and not in any finite neighbourhood thereof.

What is true in two dimensions is that, under rather broad conditions, it is possible to choose a coordinate system $(\xi,\eta)$ over a finite patch containing $P$ in which  the metric has the form
\begin{equation}\label{eq1}
(ds)^2=\lambda(\xi,\eta)^2\big((d\xi)^2+(d\eta)^2\big)\,.
\end{equation}
Locally, this is equivalent to a rescaling of Cartesian coordinates $x_1=\lambda\xi$, $x_2=\lambda\eta$, 
but where the scale factor $\lambda(\xi,\eta)>0$ is position dependent. Such a coordinate system is called \em isothermal\em.\footnote{This result and its many proofs, all quite technical, have a long history, beginning with Gauss (1822).} 

This result, by itself, has no immediate implications for our problem. A coordinate change is merely a relabelling of the points of the manifold, so it does not change measurable properties of the manifold. It preserves not only angles but also distances and areas, as long as they are correctly expressed in the new coordinates. However,
we may equally view $(\xi,\eta)$ as Cartesian coordinates on a different, planar, manifold ${\mathbb R}^2$, endowed with a Euclidean metric $(ds)^2=(d\xi)^2+(d\eta)^2$. Associating the point $(\xi,\eta)$ in this plane with the corresponding point on the curved manifold $C$ thus gives a locally invertible mapping between the two. In  general, multiplying the metric by a non-uniform scalar $\lambda(\xi,\eta)^2$ is called a \em Weyl transformation\em. It is also \em conformal\em: locally it corresponds to a rotation combined with a scale transformation, and so it preserves relative \em angles\em. Two curves in ${\mathbb R}^2$ which intersect at a point subtend the same angle between their tangent vectors at that point as do their images in the curved manifold $C$. This extends to any simply connected region $D\subset{\mathbb R}^2$  (in which all closed loops may be continuously shrunk to a point) and its image in $C$. However a Weyl transformation does not, in general, preserve lengths or areas. The combination of passing to an isothermal  coordinate system in $C$ with scale factor $\lambda$, followed by a Weyl transformation which cancels this out, gives a \em conformal mapping \em $D\to C$, which is, at least locally, invertible. As long as $D$ and $C$ have the same topology (as is the case for Hawksmoor's ceiling and a rectangle, since both are simply connected), this mapping may be extended to the whole of their interiors.

If, therefore, we take any simply connected region $D$ which is tileable by  identical squares, their corners will  map conformally to the vertices of a complete coffering of the image in $C$. The coffer lines are the images of the edges of the tiles, and therefore intersect orthogonally, and the boundary of $C$ is a union of coffer edges. Moreover the coffers are approximately square (the deviation is discussed in Sec.~\ref{secmath}). Note that the coffer lines do not have to correspond to $\xi$ or $\eta$ being constant: their pre-image in A can be any rotation of the coordinate axes. However, this does not in general give a complete coffering: there are gaps along the boundary.

Therefore, to solve Hawksmoor's problem for a given ceiling shape, we should find an appropriate set of isothermal coordinates. The mathematical results cited above imply that this is \em always \em
possible, with the coffer lines chosen to be curves of constant $\xi$ and $\eta$. The linear size of each coffer is then proportional to $\lambda(\xi,\eta)$. In fact a corollary of the Riemann mapping theorem is that there is an infinity of conformal mappings of the interior of a simply connected region $D$ of the plane to itself, and therefore the conformal mapping from this to a simply connected region of $C$ is not unique until suitable boundary conditions are imposed.  In practice, finding just one set of isothermal coordinates is a difficult, although solvable,  problem, but, as we shall see shortly,
it simplifies considerably in cases where the metric is independent of one of the coordinates, which is the case when the ceiling is a surface of revolution.

However, staying with the more general case for the time being, almost all of these conformal mappings become singular at some point(s) on the boundary, where $\lambda$ vanishes and the coffer size shrinks to zero, so there must be an infinite number of coffers. If there is to be only a finite number, the conformality of the mapping must also extend to include the boundary. Thus the boundary  of the coffered portion of the ceiling must consist of sections which are finite unions of coffer edges along which either $\xi$ or $\eta$ is constant, meeting at interior angles of $90^\circ$ or $270^\circ$.
Fortunately, as long as we ignore the truncations due to the window and door arches, the Buttery ceiling satisfies this condition, being the conformal image of a rectangle, with the N-S arches corresponding to its upper and lower edges, and intersecting the E-W sections of the horizontal cornice, which are images of the sides of the rectangle, at $90^\circ$. However for the coffering to cover the whole ceiling, the aspect ratio of the length to the width of the rectangle must be a rational number $M/N$, giving only discrete possibilities for the ratio of dimensions of the base of the ceiling.

\subsection{Immersion metric.}\label{secimm}

So far we used only results on the \em intrinsic \em geometry of 2-manifolds. However the Buttery ceiling is immersed in the 3-dimensional space of All Souls College, presumably Euclidean, and therefore inherits a notion of metric from this immersion. To specify this we may choose a Cartesian coordinate system $(x,y,z)$ such that its base, the oval boundary plane in which it meets the cornice, corresponds to $y=0$, so that $y$ measures the height of the ceiling above this plane at a particular point $(x,z)$. For an oval base which is symmetric under N-S and E-W reflections it is convenient to choose the $x$ and $z$ axes to lie parallel to the minor and major axes of the oval respectively (see Fig.~\ref{cylcoord}). 
It is more convenient, however, to use cylindrical polar coordinates $(\rho,\phi,z)$ with $\rho\geq0$, $0\leq\phi\leq\pi$ such that
$x=\rho\cos\phi$, $y=\rho\sin\phi$ and $\rho=\sqrt{x^2+y^2}$. The equation defining the immersion of the ceiling before the decoration is added then has the form $\rho=\rho(z,\phi)$, assumed to be suitably smooth. 

The three-dimensional Euclidean metric  in cylindrical polar coordinates is 
\begin{equation}
ds^2=dx^2+dy^2+dz^2=dz^2+d\rho^2+\rho^2 d\phi^2\,.
\end{equation}
Restricting this to the 2-manifold $\rho=\rho(z,\phi)$ then gives 
\begin{equation}
ds^2=dz^2+\rho(z,\phi)^2d\phi^2+\big((\partial_z\rho)dz+(\partial_\phi\rho)d\phi\big)^2\,.
\end{equation}
In general, bringing this to the isothermal form (\ref{eq1}) is difficult unless either 
$\partial_z\rho$ or $\partial_\phi\rho$ vanish. Choosing the latter means assuming that the ceiling is a surface of revolution, generated by rotating the curve $\rho=\rho(z)$ about the polar axis. The Buttery ceiling appears to satisfy this condition (small violations of which may be treated perturbatively.)
Then we can write
\begin{equation}
ds^2=\rho(z)^2\left(\frac{(1+\rho'(z)^2)}{\rho(z)^2}\,dz^2+d\phi^2\right)\,.
\end{equation}
Comparing with (\ref{eq1}) we see that a possible pair of isothermal coordinates is $(\xi,\eta)$, where $\eta=\phi$ and $\xi(z)$ satisfies
\begin{equation}\label{eq2}
d\xi/dz=\frac{\sqrt{1+\rho'(z)^2}}{\rho(z)}\,,\quad{\rm so\ that}\quad
\xi=\int_0^z\frac{\sqrt{1+\rho'(z')^2}}{\rho(z')}dz'+{\rm const.}
\end{equation}

This means that an $M\times N$ tessellation of a rectangle with identical squares will be mapped to a coffering of the ceiling in which all the coffers along an E-W row $z=$ constant are congruent, differing only by an azimuthal rotation and forming an arch. The azimuthal coordinate of each vertex takes the values $\phi=\pi n/N$, with $0\leq n\leq N$, while $\xi=\pi m/N$ with $0\leq m\leq M$.
The scale factor determining the relative size of each coffer is $\rho(z)$, the height of the  ceiling along its centre line $\phi=\pi/2$, as measured above its base. The other meridian lines are all images of this line under rotation. 
Note that, instead of imposing the  condition that the boundary is a union of coffer lines, we could extend the entire geometry by adjoining the reflection of the ceiling in the base plane $y=0$. The manifold is now a full surface of revolution, invariant under  azimuthal rotations, and the original boundary condition is equivalent to the requirement that the tiling be regular across $y=0$. Rotational invariance then implies that the solution may be expanded in terms of the sine and cosine of multiples of $2\pi/2N=\pi/N$, giving the same expressions as before.

The meaning of Eq.~(\ref{eq2}) becomes evident when expressed in terms of the infinitesimal line element along this curve $dl\equiv\sqrt{1+\rho'^2}\,dz$:
\begin{equation}
\rho\,d\xi=dl\,.
\end{equation}
Thus, if $\xi$ is a cooordinate labeling a regular square planar grid, so that $\xi=(\pi/N)\times$ integer,  the N-S length of any coffer $\int_{\rm coffer}dl$ is $\pi/N\times$ its mean height $\int\rho d\xi$ above the base, almost the same as its mean E-W width
$(\pi/2N)(\rho_m+\rho_{m+1})$. In fact these remain equal if $\rho$ varies linearly with $\xi$. Thus any violation must arise at the level of the second derivative. In isothermal coordinates this gives the local Gaussian curvature $\kappa$, with dimension (length)$^{-2}$. This to be compared to the typical area $(\pi R/N)^2$ of a coffer. Thus the relative deviations from squareness are $O(\pi^2\kappa R^2/N^2)$. For a sphere $\kappa R^2=2$, so the relative deviation is $\sim \pi^2N^{-2}$.

\subsubsection{Asymptotic behaviour.}
It is interesting to analyse the behaviour near the N or S end of the ceiling, as, say, $z\to R$. If $\rho(z)$ vanishes linearly there at an angle $\alpha$, say $\rho\sim (\tan\alpha)(R-z)$, then $\xi\sim-(1/\sin\alpha)\ln(R-z)$, so $R-z_m\propto  e^{-\pi m\sin\alpha/N}$ as $m\to\infty$.  Also if $\rho'(z)\to\infty$ as $z\to R$ (as for a hemisphere) then
$\rho_m\propto  e^{-\pi m/N}$. This means that the sequence $\{\rho_m\}$ is in geometric progression, as has been speculated. \cite{Fle2003} However according to the conformal hypothesis, this should only be the case either for a perfect cone, or for large enough values of $m$ which may not be accessible in practice. Moreover the constant or proportionality is the $N$th root of $e^{-\pi}$, not a simple rational or Pythagorean number as suggested.\footnote{This transcendental number appears naturally in physical applications of conformal symmetry. See, for example \cite{Car1984}.}

In fact this gives a strict upper bound on $\rho_m$. We have, from Eq.~\ref{eq2},
\begin{equation}
\xi<\int_0^z\frac{|\rho'(z')|}{\rho(z')}dz'=-\ln\big(\rho(z(\xi))/\rho(0)\big)\,,
\end{equation}
so that $\rho(z(\xi))\,e^\xi$ is bounded by its value at $\xi=z=0$, and therefore $\rho_m$ must decay faster than $e^{-\pi m/N}$, 
otherwise the coffering cannot be conformal.

\subsection{Examples.}
\subsubsection{Cone.}\label{seccone}
As discussed above, this leads to the sequence $\{\rho_m\}$ being in geometric progression, with ratio $e^{-\pi(\sin\alpha)/N}$, where $2\alpha$ is the opening angle.
\subsubsection{Hemispherical ceiling.}\label{sechemi}
Suppose the ceiling is hemispherical, so that its base is a circle of radius $R$. In cylindrical polar coordinates it is given by $\rho=\sqrt{R^2-z^2}$ with $-R\leq z\leq R$ and $0\leq\phi\leq\pi$. Eq.~\ref{eq2} becomes
\begin{equation}\label{hemi}
d\xi/dz=\frac R{R^2-z^2}\\,\quad\text{so that} \quad
\xi=-(1/2)\ln\left(\frac{R-z}{R+z}\right)\,.
\end{equation}
Solving for $z(\xi)$,
\begin{equation}\label{zxi}
z(\xi)=R\tanh\xi\,,
\end{equation}
and the scale factor is $\sqrt{R^2-z^2}$ so that the coffers become smaller in the polar regions.

Taking a regular $M\times N$ square array in the $(\xi,\phi)$ with $\xi=m\pi/N$, $\phi=n\pi/N$ where $-M/2\leq m\leq M/2$, $0\leq n\leq N$, its 3d image on the ceiling  is
\begin{equation}\label{ceilingvertices}
\{\phi=\pi n/N, z=z_m=R\tanh(\pi m/N),\rho=\rho_m=R\,{\text{sech}}(\pi m/N)\}\,.
\end{equation}
This is shown in Fig.~\ref{recon}. Note that in this example there is a simple iterative procedure for generating the $\{z_m\}$, since
\begin{align}
z_{m+1}&=R\tanh\big(\pi(m+1)/N\big)
=R\frac{\tanh(\pi m/N)+\tanh(\pi/N)}{1+\tanh(\pi m/N)\tanh(\pi/N)}\,,\\
&=\frac{z_1+z_m}{1+z_1z_m/R^2}\,,  \label{iter}
\end{align}
where $z_1=R\tanh(\pi/N)\approx \pi/N\ll 1$. In fact more generally
\begin{equation}
z_{m+n}=\frac{z_m+z_n}{1+z_mz_n/R^2}\,.
\end{equation}

Writing $z=R\sin\theta$, (\ref{hemi}) is one of several equivalent expressions\cite{Sny1993} for Mercator's projection of a sphere, with $\theta$ being the latitude in radians, so that, as claimed, Hawksmoor's coffering is its inverse.

\subsubsection{The Buttery ceiling}\label{secbutt}
This is only slightly more complicated. We have
\begin{align}\label{butt}
\begin{split}
\rho(z)&=\sqrt{R^2-(z-R')^2}\,\qquad(z\geq R')\\
&=R\,\qquad(-R'\leq z\leq R')\\
&=\sqrt{R^2-(z+R')^2}\,\qquad(z\leq -R')\,,
\end{split}
\end{align}
where the length of the barrel section of the ceiling is $2R'$. If there are $N$ coffers in each row, looking at the hemispherical section of the ceiling we see that each central coffer has width, and therefore also length, $\approx\pi R/N$. $2R'$ must be an integer multiple of this. Therefore $2R'/R$ must be a rational multiple of $\pi$. In  actuality it is $4\pi/11$. So the rectangular base of the barrel section has aspect ratio $R'/R=2\pi/11\approx 0.57$, which agrees with observation.
The solution is
\begin{align}\label{hc}
\begin{split}
z_m&=R\tanh((m-m')\pi/N)+Rm'(\pi/N)\,,\quad(m\geq m')\\
&=Rm(\pi/N)\,,\quad(-m'\leq m\leq m')\\
&=R\tanh((m+m')\pi/N)-Rm'(\pi/N)\,,\quad(m\leq -m')\,,
\end{split}
\end{align}
where $N=11$ and $m'=2$. The actual coffering according to this prediction is illustrated in Fig.~\ref{recon}.

In order to compare with the actual ceiling in the absence of direct measurements, we should simulate the effect of the camera projection on the predicted 3d data, and compare the result with actual photographs. Fortunately this is straightforward, as, in the limit when it is much smaller than the object, it may be modelled by a simple pinhole camera, as shown in Fig.~\ref{cylcoord}. At this level all that a lens does is to focus different rays from a particular point on the ceiling to the same point in the focal plane, thus allowing for a wider aperture, as well as to magnify or diminish the whole image by a uniform scale factor, which is not important as we are interested only in ratios of lengths.

We assume, as appears to be the case, that in Fig.~\ref{fig4} the camera is located directly below the apex of the ceiling, and directed vertically upwards.
After some simple geometry we find that if the Cartesian coordinates of the object $P$ are $(x,y,z)$ (see  Fig.~\ref{cylcoord} for notation), those of its image $P'$ are
\begin{equation}
x'=-\frac{fx}{h+z}\,,\qquad y'=-h-f\,,\qquad z'=-\frac{fz}{h+z}\,,
\end{equation}
from which we see that the only relevant parameter is the ratio of $h$ to the height of the ceiling above its base. 

The result of applying this transform, overlaying the actual photograph, is shown in Fig.~\ref{comp}, in which $h/R$ has been set to 1.2 for the best overall fit, corresponding to the camera having been positioned on or near the floor.

\subsubsection{The Pantheon dome.}\label{secpan2}

The analysis for the Pantheon dome, assumed to be hemispherical, is almost identical, except that the $z$-axis of symmetry now lies vertically. Also, because the base of the dome is a circle,
the mapping is now from a finite cylinder, with $\pi/N$ replaced by $2\pi/N$ where $N$ is now the number of coffers in each $360^\circ$ row ($N=28$ in this case). Thus the vertices lie at 
\begin{equation}
\phi=2\pi n/N\,, \qquad z=z_m=R\tanh(2\pi m/N)\,,\qquad \rho=\rho_m=R\,{\text{sech}}(2\pi m/N)\}\,,
\end{equation}
while the scale factor is $\propto {\text{sech}}(2\pi m/N)$.
In this case there are direct measurements, reported in \cite{Ali2017}, of the 
row dependence of the mean widths and heights of the coffers (as measured along the surface,) which should, according to the conformal hypothesis, both be proportional to this scale factor.  Thus their ratios should be universal. The comparison is shown in Fig.~\ref{ratios}. This appears encouraging, since there is no adjustable parameter. However, with only four values of $m$ for comparison, which do not explore the exponential tail of the distribution, one may only assert that these results do not refute the conformal hypothesis.

Stronger confirmation comes from a comparison with photographs, as for the Hawksmoor ceiling. In this case, assuming that the camera is directed vertically upwards, directly below the centre of the oculus, a study of the geometry similar to that in Fig.~\ref{cylcoord} shows that the polar coordinates of the camera image of the vertices are $\phi'=\phi+\pi$ and
\begin{equation}
\rho'_m=\frac{f\rho_m}{z_m+h}
=\frac{fR\,{\text{sech}}(2\pi m/N)}{R\tanh(2\pi m/N)+h}
\propto \frac1{\sinh(2\pi(m+m_0)/N)}\,,
\end{equation}
where $h/R=\tanh(2\pi m_0/N)$. The comparison is shown in Fig.~\ref{pancompph2}.
Apart from the overall scale, there is just one parameter $m_0$, which was chosen to obtain the best overall fit. It corresponds to $h/R\approx 0.95$. Given that the height of the base of the dome above the floor is stated to be one half its diameter of 43.3 m, this corresponds to the camera being positioned at a height of approximately one meter. 

\section{Why the Pantheon coffering is conformal.}\label{secDir} 
\gdef\p{{\partial}}

The conformal hypothesis, described in the previous sections, is an economical way of accounting for the observed cofferings of both the Pantheon and the Hawksmoor ceilings, as well as other examples. However, it is a somewhat abstract mathematical condition, and it would be more compelling if it were to follow from some  general physical principle. We now argue that this is so. This section is rather more mathematical than the earlier ones.

Let us assume that we know nothing about the construction methods of the Romans except that they were not gods -- they were bound by the laws of physics, even though they did not possess our modern understanding of them. Instead, they were guided by practical rules which, experience showed, produced stable structures.

The modern view of the laws of classical physics is that they are encapsulated in variational principles, which identify functionals of all the possible degrees of freedom and their corresponding velocities, whose first order variation vanishes whenever the laws are being followed. For the case of static equilibrium this is the \em principle of virtual work, \em which states that if a physical system is in equilibrium, the total work done by all external forces and moments is zero for any arbitrary, infinitesimally small hypothetical displacement. In the present context, this means that the total potential energy, which is the sum of the elastic energy and gravitational potential energy, should be extremized over all possible deformations. Those   virtual deformations which do not extremize this energy do not obey the laws of classical mechanics. In the context of coffering a ceiling, varying the map $f:D\to C$, that is the coffering design, is a deformation. 

Another important notion in the physics of continuous media is the idea of an effective energy functional, which governs the long wavelength coarse-grained degrees of freedom, or soft modes, which are often the most critical for stability.  It has form of an integral over space of an energy density, which is given by keeping the first nontrivial terms in a gradient expansion, usually quadratic in the first derivatives.  The form of these is often fixed by global symmetries, and depends on the microscopic details through only a few phenomenological coefficients. In that spirit let us construct an effective energy governing the soft modes of the coffering design. 

Suppose, as described in Sec.~\ref{secmain}, the ceiling is parametrized by coordinates $x=(x^1,x^2)\in D\subset{\mathbb R}^2$ through the smooth bijection 
$f: x\leftrightarrow \vec X(x)\in {\mathbb R}^3\cap C$.  Each such mapping then gives rise to a (sequence of) coffering design(s) of $C$ with curvilinear quadrilaterals. It is a Hawksmoor coffering if $f$ is conformal. 

Let $T_XC$ be the tangent plane at a point $X$ on $C$. It is orthogonal to the normal vector $\vec n_X$, with respect to which $\vec X(x)$ may be decomposed into a normal component $X_n=\vec n.\vec X$ and a projection $\vec{\widetilde X}(x)=\vec X(x)-X_n(x)\vec n\in T_XC$. Thus $(x^1,x^2)$ gives a smooth parametrisation of a patch of $T_XC$ around the point where it touches $C$. On the other hand, one may choose a Cartesian basis for ${\mathbb R}^3$ at $X$ such that $X_3=X_n$ and $\widetilde X^\gamma$ with $\gamma=1,2$ are Cartesian coordinates on $T_XC$.\footnote{In this section the summation convention is assumed, with Greek letters running from 1 to 2 and Latin letters from 1  to 2 or 3, depending on the context.  We also get serious in distinguishing covariant and contravariant vectors by upper and lower indices, with Greek indices lowered by $g$ and Latin indices by the Cartesian metric on ${\mathbb R}^3$, more generally some 3-dimensional Riemannian manifold.}

Since, by the chain rule, 
\begin{equation}
ds^2=dX^idX_i=(\p\widetilde X^i/\p x^\alpha)(\p\widetilde X_i/\p x^\beta)dx^\alpha dx^\beta +dX_n^2\,,
\end{equation}
the immersion induces a metric $g_{\alpha\beta}=
(\p\widetilde X^i/\p x^\alpha)(\p\widetilde X_i/\p x^\beta)$ on $T_XC$. Note that as $X$ moves on $C$, the tangent plane in general rotates. However the manipulations which follow depend on the  metric $g(X)$ only at $X$ (and not its derivatives), which then may be regarded as a \em uniform \em flat metric on $T_XC$. 

The local properties  of the coffering design are now encoded in the $2\times2$ transformation matrices $X_{i,\alpha}=\p\widetilde X_i/\p x^\beta$, and it is from these that the effective energy $E[f]=\int e(x,f)dx^1dx^2$ is built. The density $e(x,f)$ is required to depend only on $X$ and quadratically on its first derivatives $X_{i,\alpha}$, and be invariant under the various symmetries of the problem. Any terms which involve higher powers of the first derivatives, or higher order derivatives, are irrelevant for the long wavelength 
modes.

Let us write out explicitly the various possible independent terms in $e$, in descending order of their degree of symmetry:
\begin{eqnarray}
e_D=|X_{1,1}|^2+|X_{1,2}|^2+|X_{2,1}|^2+|X_{2,2}|^2\,, \label{eDir}\\
(X_{1,2}-X_{2,1})^2\,, \label{erot}\\
(X^1_{,1}+X^2_{,2})^2\,,\label{ecomp}\\
X_{1,1}X_{2,2}-X_{1,2}X_{2,1}\,,\label{edet}\\
X_{1,1}^2+X_{2,2}^2, X_{1,1}X_{2,2}\,,\label{ecub}
\end{eqnarray}
where $|X_{i,\alpha}|^2$ is shorthand for $X_{i,\alpha}X^{i,\alpha}$ (not summed). 
Each of these, on integration $\int dx^1dx^2$, is invariant under scale transformations $x^\alpha\to bx^\alpha$ for $b>0$. In addition, 
the first term (\ref{eDir}) is known as the Dirichlet energy density, and is invariant under separate O$(2)$ rotations in $X$ and $x$.\footnote{In fact it is conformally invariant, as may be seen by writing the integral in complex coordinates $(z,\bar z)=(x^1\pm ix^2)$, when it becomes $\int (g_{z\bar z}/\sqrt{g_{zz}g_{\bar z\bar z}})dA$, which is invariant under any complex analytic reparametrizations $(z,\bar z)\to (w(z),\bar w(\bar z)$. Thus $E_D$ depends only on the conformal equivalence class of the map $f$.} (\ref{erot}) and (\ref{ecomp}), respectively, measure the relative local rotation and dilatation between $X_i$ and $x_\alpha$,  and are invariant under simultaneous O$(2)$ rotations in both $X$ and $x$. The last two terms in (\ref{ecub}) break O$(2)$, and retain only the ${\mathbb Z}_4$ symmetry of the square lattice, but are essential because they are the only contributions at this level of description which couple the long wavelength deformations to the orientation of the square coffering. Depending on the sign of the coupling constant, this will favour  the coffers aligning along the isolines of $(x^1,x^2)$, or being relatively rotated by 45$^\circ$.

The above expressions are familiar from linear elasticity theory\cite{LL}, with $X$ replaced by the displacement field $u$, and in the limit when $C$ becomes flat, they are identical up to a total derivative. It is tempting, therefore, to identify the above densities as contributing to the physical elastic energy of the structure, the first three representing shear, vortex, and compressional energies. However their significance here is that minimising each of them leads to the map
$f$ necessarily being conformal. This is a well known result\cite{Dir} for the Dirichlet energy, but it holds more generally.

Reorganising the above expressions as
\begin{eqnarray}
(X_{1,1}-X_{2,2})^2+(X_{2,1}+X_{1,2})^2
+2\det(\p X/\p x)\,,\label{eDir2}\\
(X_{1,2}-X_{2,1})^2\,, \label{erot2}\\
(X_{1,1}-X_{2,2})^2+4\det(\p X/\p x)+4X_{2,1}X_{1,2}\,,\label{ecomp2}\\
X_{1,1}X_{2,2}-X_{1,2}X_{2,1}=\det(\p X/\p x)\,,\label{edet2}\\
(X_{1,1}-X_{2,2})^2+2\det(\p X/\p x)+2X_{2,1}X_{1,2}\,.\label{ecub2}
\end{eqnarray}
Any linear combination of these with non-negative coefficients is bounded from below by a non-negative multiple of the Jacobian $\det(\p X/\p x)$ which expresses the area element $dA$ of $C$ in terms of the measure $dx^1dx^2$. Moreover, this is achieved when $X_{1,1}=X_{2,2}$ and $X_{1,2}=X_{2,1}=0$. In terms of the induced metric $g$, this implies that $g_{11}=g_{22}$ and $g_{12}=0$, that is, the coordinates $x^\alpha$ are isothermal and the map $f: D\to C$ is conformal. 

The normal component $X_n=\vec X.\vec n$ has so far been overlooked. It gives a contribution $\propto(\p_\alpha X_n)(\p^\alpha X_n)$ to the energy density, where now the derivatives act on $\vec n$, and so involve the curvature of $C$. In fact it is proportional to the squared mean curvature
\begin{equation}
H^2=(1/4)(\kappa_1+\kappa_2)^2=(1/4)(\kappa_1-\kappa_2)^2+\kappa_1\kappa_2\,,
\end{equation}
which may be identified with the bending energy density.
For a conformal metric  $\kappa_1=\kappa_2$, and so this is bounded from below by the Gaussian curvature $K=\kappa_1\kappa_2$. However the integral of this over a closed manifold, or an open one with appropriate boundary conditions, is proportional to the Euler character of the manifold, and is therefore independent of the map $f$.

Thus conformality, in this context, is an emergent property arising from applying the principle of virtual work to the elastic energy.
An important feature of the above inequalities is that they hold point-wise, so the integrands can be multiplied by the same non-negative function of $X$ (but not its derivatives) and the corresponding inequality is still saturated if the map is conformal. 

The gravitational potential energy is
\begin{equation}\label{egrav}
E^{gr}[X]=-\int\rho(X)\vec g\cdot(\vec X)dA=-M\vec g.\vec X_c\,,
\end{equation}
where $\vec g$ is the magnitude and direction of the acceleration due to gravity, $\rho$ is the mass per unit area, and 
$\vec X_c$ is the location of the centre  of mass.  For a fixed a total mass $M$, $E^{gr}$ is independent of the map $f$. The architects are instructed to make plans for constructing a fixed ceiling $C$ with a certain mass of material. They also have to obey the laws of physics, and therefore should minimize $E^{el}[X]$ over all smooth maps $\hat f:D\to C$, subject to the constraints that $\hat f(D)=C$ and the total mass $M$ is fixed. In a variational approach the constraint may be implemented by adding a term $-\Lambda\int\rho(X)dA$ to the elastic energy, with a Lagrange multiplier $\Lambda$. This is just the gravitational energy in (\ref{egrav}), identifying $\vec G.\vec X_c$ with $\Lambda$. The Euler-Lagrange equations for varying $X(x)$ with a fixed parametrisation are then linear, and imply that each $X_i(x)$ is a harmonic function. However one should actually vary the parametrisation for fixed $X\in C$. The two are related in the standard way by a Legendre transform. 

However, there is a finite degeneracy of solutions, due to the  residual ${\mathbb Z}_4$ symmetry. This is broken by the boundary conditions. In the example of the Buttery ceiling, we argued that it is then unique. Finally, the gravitational energy is minimized by lowering the centre of mass as much as possible. This may be achieved by having the thickness of the shell decrease with increasing height, as actually occurs in the Pantheon dome\cite{Wad2008}, although it should not be too thin or the local stress will become too large. 

To summarise this section, the principle of virtual work implies that any equilibrium must correspond to an extremum of the total potential energy, which should be a local minimum for stability against finite deformations. We have argued that the contribution to this from long wavelength elastic shear and compression modes is given by a generalisation of standard linear elasticity, which is minimised when the map from the plane to the ceiling is conformal. The boundary conditions then make this unique. A coffering design based on this map would be the only one completely stable to long wavelength shear and bending deformations. Any other scheme would correspond to an extremum with at least one unstable mode. Therefore one might say that if the coffering had not been close to conformal, the structure would not have lasted nearly 2,000 years.

\section{How the Romans might have done it.}\label{romhow}

The builders of the Pantheon were probably unaware of the principle of virtual work, although the concept was known to Hero of Alexandria (1st-2nd c. CE), and the foundations of statics had been established by Archimedes (287-212 BCE). No doubt they developed their own sophisticated recipes for construction with concrete. Nevertheless they were bound by the laws of physics, and, possibly by trial and error, converged on a conformal coffering scheme in order to balance the forces, which, as argued in the previous section, is a consequence of minimizing the sum of the elastic and gravitational potential energies while keeping the total mass fixed. Experience may have taught them that a conformal coffering with (approximate) squares,
whose edges lie along meridian and parallel lines, gives the most stable structure, as argued in the previous section. 

Alternatively, they may have been drawn to this scheme for aesthetic or other reasons.  In any case, they would be trying to maximize the `squareness' of each coffer. One way might have been to estimate this is through a discretized version of the Dirichlet energy\cite{Des2002}. This is an important problem in computer science and animation, identifying the best way to digitize a 2d surface immersed in three dimensions, although triangulations are usually used in preference over quadrangles.\cite{Pink1993}. Consider a quadrangulation of the ceiling in which the vertices are labelled by integers $(m,n)$. Each vertex is then assigned a location $\vec X_{m,n}$ on the ceiling, ordered so that nearest neighbours to this vertex carry labels $(m\pm 1, n\pm 1)$.  This gives a coffering of the ceiling. A discrete approximation to the Dirichlet energy is
\begin{equation}
E_D^{dis}[X]=(1/2)\sum_{m,n}\left[|\vec X_{m+1,n}-\vec X_{m,n}|^2+|\vec X_{m,n+1}-\vec X_{m,n}|^2\right]\,,
\end{equation}
or, simply put, the sum of the squared lengths (as measured in 3d) of all the nearest neighbour edges. It is easy to see that the summand is bounded from below by the area. Writing it as $|\vec P|^2+|\vec Q|^2$, 
\begin{equation}
(|\vec P|^2+|\vec Q|^2)^2=(|\vec P|^2-|\vec Q|^2)^2+4|\vec P|^2|\vec Q|^2\geq
4|\vec P.\vec Q|^2+4|\vec P\times\vec Q|^2\geq 4|\vec P\times\vec Q|^2\,,
\end{equation}
where the last term is the square of the area of the triangle with vertices at $(X_{m,n},X_{m+1,n},X_{m,n+1})$. Moreover, the bound is attained when this is both isosceles  and right-angled. Thus the ratio of the area of each coffer to its mean square edge length gives a measure of its `squareness', the maximum value of 1 being attained for a perfect square immersed in ${\mathbb R}^3$. 

Minimizing $E_D^{dis}$ subject to the constraint that all vertices lie on $C$ is straightforward with modern computers, and could have been carried out approximately, adding one coffer at a time, by the Romans. Indeed, this appears to be the prescription suggested by Aliberti and Alonso-Rodr\'iguez \cite{Ali2017}, quoted in Sec.~\ref{secrec}.

However it is also the case, for a conformal coffering, that an observer looking directly along a normal to the surface will see the coffer directly in front as a slightly deformed square, with curved edges that however are all of equal length (as measured in three dimensions), and which meet at right angles. Thus an observer positioned at the center of the circular base  could, in principle, have communicated to the builders
when each new coffer looked as square as possible.

\section{How Hawksmoor might have done it.}\label{sechow}
Of course, Hawksmoor probably did not know about hyperbolic functions even as used earlier by Mercator, and certainly nothing of differential geometry. Equally it is difficult to guess the iterative solution (\ref{iter}), so it is interesting to speculate how he might have constructed his coffering to such a high degree of accuracy. The first point to realise is that, unlike most other coffers made of several different materials including wood, the Buttery ceiling is composed of individual pieces of masonry, one for each coffer square, carefully sculpted on the lower, visible face and the sides, but left uncut on its upper face. Probably Hawksmoor realised that all the coffers in each row should be identical in size and shape, and in fact the coffers in all four of the central rows are all the same. He could therefore have his masons mass produce these. Assuming the window and doorway arches to be in place, the four central  arches of coffers could be erected in a manner no doubt very familiar to builders of that period. Where the ceiling becomes a quarter sphere, it would have been necessary systematically to reduce the size of each coffer and to taper them slightly. If Hawksmoor had realised that their size should diminish proportional the height of the apex of each arch above that of the cornice, this would have been useful. This protocol is more or less equivalent to solving (\ref{eq2}) by discretization. Such approximations can lead to instabilities  in the numerical solution. However, as we have pointed out, a coffering is possible for any shaped ceiling which is surface of revolution, so it would not have been necessary to achieve great accuracy, as long as the coffers within each arch were sufficiently similar. We note that this protocol is much simpler than some of those suggested for the Pantheon dome, largely because, if it were designed correctly, gravity would stabilise each arch as it was completed. This would not have been the case for the Pantheon dome, whose construction was a masterpiece of engineering and design for the period.

\section{Conclusion.}\label{secconcl}
The patterns of coffers on both the Hawksmoor ceiling in All Souls and the Pantheon dome appear to be consistent with the hypothesis that they are conformal images of a regular square tiling of a flat surface. Indeed, we have argued that this solution  guarantees that the coffer lines, as well as their diagonals, intersect orthogonally, and that, if the coffering is complete, the conformal mapping is unique, although the explicit results depend on the ratio $M/N$ of the numbers of coffers in each row and column. It also implies that the mapping from the flat to the curved surface is  the mathematical inverse of Mercator's canonical projection. 

It should be stressed that, by itself, this is only a mathematical hypothesis which accounts for the most obvious features of the coffering, that the coffer ridges intersect orthogonally and the coffers themselves are almost square. It has the advantage of simplicity as compared with other analyses, and, all other things being equal, should therefore be preferred on the basis of Occam's razor. Against the obvious objection that Hawksmoor, never mind the Romans, knew nothing of differential geometry, we have argued that the conformal mapping approach leads to a unique solution which may well be accessible by more elementary methods, as suggested in Secs.~(\ref{romhow},\ref{sechow}). 

However, in Sec.~\ref{secDir} it was shown that  conformality is a generic consequence of minimising a candidate expression for the continuum elastic energy of the soft modes, as well as a discrete version  of the same. If the coffering design were not conformal, the structure would be unstable to some long wavelength modes. 

Meanwhile, may the Warden and Fellows of All souls continue to enjoy their repast, in the knowledge that their Buttery ceiling, being conformal, still stands proudly in a lineage of pure and applied scientific knowledge extending from Hadrian's Empire, through the golden age of global exploration and the charms of the Baroque, to the revolutionary mathematics of the 19th century which, to this day, continues to inspire and inform the work of its Fellows.

\enlargethispage{20pt}  

\noindent\em Acknowledgements\em. The author is  grateful to Alison Turnbull for rekindling his interest in the Hawksmoor ceiling and to Anthony Geraghty for sharing his expert knowledge of Hawksmoor and his All Souls projects, as well as suggesting the comparison with the Roman Pantheon, He thanks Peregrine Horden and C\'ecile Fabre for encouraging this project, Steve Evans and Afshan Sohail for providing survey materials. Some of the figures in this paper were produced using Mathematica\texttrademark.

\end{document}